\documentclass[11pt]{amsart}
\usepackage{amsmath}
\usepackage{amssymb}
\usepackage{tabularx}
\usepackage{enumerate}
\usepackage{graphicx}
\usepackage{texdraw}
\usepackage{color}

\usepackage{mathrsfs}

\usepackage{amsfonts,amssymb,amsmath}
\usepackage{epsfig}

\makeatletter
\@addtoreset{equation}{section}

\makeatother
\newcommand{\R}{\mathbb{R}}

\begin{document}
\title[Singular Limits of Porous Media Equations with Reactions]
      {Singular Limits of Porous Media Equations with Bistable Reactions$^*$}
\thanks{$^*$ This research was partly supported by Natural Science Foundation of China (No. 12071299).}
\author[B. Lou]{Bendong Lou$^{\dag}$}
\thanks{$\dag$ Mathematics and Science College, Shanghai Normal University, Shanghai 200234, China.}
\thanks{{\bf Emails:}  {\sf lou@shnu.edu.cn (B. Lou)}}
\date{}


\maketitle

\begin{abstract}
We consider a porous media equation with balanced bistable reactions, equipped with some general nonlinear boundary condition.
When the coefficient of the reaction term is much larger than that of the diffusion term, we see that, besides the possible free boundary, sharp interfaces appear between two stable steady states. By using the method of matched asymptotic expansions, we derive the motion law of each interface, which is a mean curvature flow (may depends on normal direction of the interface). In addition, the original boundary condition reduces to Robin ones at the points where the interface contacts the domain boundary.
\end{abstract}

%%%%%%%%%%%%%%%%%%%%%%%%%%%%%%%%%%%%%
%%%%%%%%%%%%%%%%%%%%%%%%%%%%%%%
%%%%%%%%%%%%%%%%%%%%%%%%%%%%%%
%%%%%%%%%%%%%%%%%%%%%%%%%%%%%%
%%%%%%%%%%%%%%%%%%%%%%%%%%%%%%

\section{Introduction}
Consider the Allen-Cahn equation
\begin{equation}\label{AC}
u_t =\Delta u + \frac{1}{\varepsilon^2} f_0 (u),
\end{equation}
where $\varepsilon>0$ is a parameter and $f_0$ is a bistable reaction term like
\begin{equation}\label{ex-f0}
f_0(u)=u \Big(u-\frac12 \Big) (1-u),
\end{equation}
which has balanced double-well potential $F_0(u):= -\int_0^u f_0(s) ds$.
It is known that, sharp interfaces (or, internal layers) often appear in \eqref{AC} when we take the singular limit as $\varepsilon\to 0$, and each interface moves according to its mean curvatures:
$$
V= - (N-1)\kappa,
$$
where $V$ is the normal velocity of the interface and $\kappa$ is the mean curvature of it (see, for example, \cite{Matano2008, Chen} for rigorous approach and \cite{Lou2007, NMHS} for formal derivation).

Assume $\Omega$ is a domain in $\R^N$ with smooth boundary. In this note we consider a nonlinear diffusion version of \eqref{AC} in $\Omega$, that is, a porous media equation (PME, for short) with a reaction term:
\begin{equation}\label{P}
\left\{
\begin{array}{ll}
\displaystyle u_t = \Delta u^m + {\bf q}(x,t)\cdot \nabla u^k + \frac{1}{\varepsilon^2}f(u),& x\in \Omega,\ t>0, \\
\displaystyle \frac{\partial u^m}{\partial \nu} = \frac{1}{\varepsilon} G(x,t,u), & x\in \partial \Omega, t>0,
\end{array}
\right.
\end{equation}
where $m\geq 1,\ k\geq 1,\ \varepsilon>0$ are parameters, ${\bf q}$ is a smooth vector function, $\nu$ denotes the unit outward normal vector on $\partial \Omega$, $G$ is a smooth function, and $f$ is a bistable reaction term satisfying:
\begin{equation}\label{cond-f}
\left\{
\begin{array}{l}
f\in C^2([0,1]),\ \ f(0)=f(a)=f(1)=0 \mbox{ for some }a\in (0,1), \\
 f(u)<0 \mbox{ for }u\in (0,a),\quad f(u)>0 \mbox{ for }u\in (a,1),\\
\int_0^1 r^{m-1}f(r) dr =0.
\end{array}
\right.
\end{equation}
We will in the next section that the last equality about the integral implies that the equation without the advection term admits a stationary solution connecting two steady states  $0$ and $1$. So it corresponds to the nonlinearity $f_0(u)$ in \eqref{ex-f0} for the Allen-Cahn equation.

When $m>1$, the diffusion term $\Delta u^m$ means that the diffusion is a density-dependent one.
Such equations can be used to model population dynamics with diffusion flux depending on the population density, the combustion, propagation of interfaces and nerve impulse propagation phenomena in porous media, as well as the propagation of intergalactic civilizations in the field of astronomy (cf. \cite{GN, GurMac, NS, SGM} etc.). The well-posedness of the problem \eqref{P} was studied in \cite{A1969, A1, Sacks, Vaz-book, Wu-book} etc.. The second term on the right hand side of the equation in \eqref{P} denotes advection, which means the reaction (nonlinear) diffusion process occur in an advective environment. The most studied case is $k=1$ or $k=m$ (see, for example, \cite{GK})

In Section 2, we will give a rough explanation on the appearance of the interface, and specify the monotone stationary solution connecting $0$ and $1$ of the equation $u_t =\Delta u^m + f(u)$. In Section 3 we will use the method of matched asymptotic expansion to show that, as $\varepsilon\to 0$, the motion law of the interface is a mean curvature flow. In Section 4, we still use the asymptotic expansion method to show that the boundary condition in \eqref{P} reduces to a Robin one for the interface.

\section{Interfaces and Monotone Stationary Solutions}
One important feature of the PMEs is the existence of free boundaries.  More precisely, assume
that $u_0(x)\in C(\R^N)$ has compact support, then the solution of $u_t =\Delta u^m$  with initial data $u(x,0)=u_0(x)$ has some free boundaries. For each point $x_0$ on the initial free boundary, it has a waiting time, and after that time the boundary moves according to the Dacry's law (cf. \cite{Vaz-book, Wu-book}).

For simplicity, from now on we consider the problem \eqref{P} with initial data
\begin{equation}\label{ini-1}
u(x,0)=u_0(x) \in C(\overline{\Omega}),
\end{equation}
with
\begin{equation}\label{ini-2}
\Omega_0^- := \{ x\in \Omega \mid 0\leq u_0(x) <a\} \not= \emptyset,\qquad
\Omega_0^+ := \{ x\in \Omega \mid a<u_0(x) \leq 1\} \not=\emptyset.
\end{equation}
Note that when $u_0(x)>0\ (x\in \Omega)$, the solution $u$ has no free boundary, otherwise, it has.
We now roughly explain that, when $\varepsilon\ll 1$, besides the possible free boundaries, the solution develops sharp interfaces.
We consider the special case where $\partial \Omega_0^+ \cap \partial \Omega_0^- \not= \emptyset$. Assume $x_0$ is one point in such a set, then there are some points $x^\pm \in \Omega_0^\pm$
near $x_0$. Since $\varepsilon>0$ is sufficiently small, the first and second terms on the right hand side of the equation in \eqref{P} can be neglected, and $u_t(x^\pm, t)\sim \frac{1}{\varepsilon^2} f(u(x^\pm, t))$.
This means that $u(x^+, t)/u(x^-,t)$ increases/decreases very fast and approaches $1/0$ in the time scale $O(\varepsilon^2)$. In such a way, $u$ develops a sharp interior layer connecting $1+o(1)$ and $o(1)$. Assume
the thickness of such a layer is $\delta$, then, roughly, the largest value of $|\nabla u^m|=O(\frac{1}{\delta})$, the largest value of $|\Delta u^m|= O(\frac{1}{\delta^2})$. In the right hand said of the equation in \eqref{P}, in order for the first and the second terms balance out the third, we require that $\delta \sim \varepsilon$. This means that, roughly, the thickness of the layer is $O(\varepsilon)$.
(In the typical Allen-Cahn equation, it was shown rigorously in \cite{Chen} that the thickness is $O(\varepsilon |\ln \varepsilon|)$.) The aim in this paper is to show, by using the method of matched asymptotic expansion, that the interface (that is, the $a$-level set of the solution) moves according to a mean curvature flow as $\varepsilon\to 0$.

Next we consider the stationary solution connecting $0$ and $1$ of a one dimensional
simplified equation, that is, consider
\begin{equation}\label{ss}
\left\{
\begin{array}{ll}
(\varphi^m)_{yy} + f(\varphi)=0, & y_* < y <y^*,\\
\varphi(y_*)=0,\quad \varphi(y^*)=1, &
\end{array}
\right.
\end{equation}
where $y_*$/$y^*$ denotes the left/right most boundary of the definition interval in which $0<\varphi(y)<1$.
According to the degeneracy of $f$ near $0,1$, $y_*$/$y^*$ may be finite number or $-\infty$/$\infty$.
Denote by $\Phi$ the unique solution of the initial value problem
\begin{equation}\label{ini-p}
\varphi'(y) = \sqrt{\frac{2}{m}} \cdot \frac{\sqrt{-F(\varphi)}}{\varphi^{m-1}},\qquad \varphi(0)=a,
\end{equation}
where
\begin{equation}\label{def-F}
F(u):= \int_0^u r^{m-1} f(r) dr,\qquad u\in [0,1].
\end{equation}
Then by the assumption in \eqref{cond-f} we see that
$$
F(u)<0 \mbox{ for }u\in (0,1),\qquad F(0)=F(1)=0.
$$
Denote the maximum existence interval of $\Phi$ by $(y_*, y^*)$, then $\Phi$ is the solution of \eqref{ss}.

To determine $y_*$ and $y^*$, we need the decay rate of $f$.  Let us just assume for a moment that $f\in C([0,1])$ in this part. Precisely, assume, for some $\alpha_0,\ \alpha_1\geq 0$,
$$
f(u)\sim -u^{\alpha_0} \mbox{ as } u\to 0^+,\qquad  f(u)\sim (1-u)^{\alpha_1} \mbox{ as } u\to 1^-.
$$
Then, as $u\to 0^+$, we have $F(u)\sim -\frac{u^{m+\alpha_0}}{m+\alpha_0}$, and so by \eqref{ini-p} we have
\begin{equation}\label{varphi-sim-eq}
\varphi'(y) \sim \sqrt{\frac{2}{m(m+\alpha_0)}} \varphi^{\frac{m+\alpha_0}{2} -m +1} \mbox{ as } y\to y_*.
\end{equation}
This implies that
\begin{equation}\label{decay-0}
\left\{
\begin{array}{l}
\mbox{when } 0< \alpha_0<m,\ \ y_*>-\infty \mbox{ and } \Phi(y)\sim O(1)(y-y_*)^{\frac{2}{m-\alpha_0}} \ (y\to (y_*)^+), \\
\mbox{when } 0<\alpha_0=m,\ \  y_* = -\infty \mbox{ and } \Phi(y)\sim O(1)e^{\frac{y}{m}} \ (y\to -\infty),\\
\mbox{when } \alpha_0 >m,\ \  y_* = -\infty \mbox{ and } \Phi(y)\sim  \left( O(1) - \frac{\alpha_0 -m}{2} y \sqrt{\frac{2}{m(m+\alpha_0)}} \right)^{\frac{2}{m-\alpha_0}}.
\end{array}
\right.
\end{equation}
In particular, $y_*>-\infty$ in the typical case where $\alpha_0=1$. Similarly, we have
$$
\left\{
\begin{array}{l}
\mbox{when } 0<\alpha_1<1,\ \ y^* < \infty \mbox{ and } \Phi(y)\sim O(1)(y^* -y )^{\frac{2}{1-\alpha_1}} \ (y\to (y^*)^-), \\
\mbox{when } \alpha_1=1,\ \  y^* = \infty \mbox{ and } \Phi(y)\sim 1- O(1)e^{- \frac{y}{\sqrt{m}}} \ (y\to \infty),\\
\mbox{when } \alpha_1 >1,\ \  y^* = \infty \mbox{ and } \Phi(y)\sim  \left( O(1) + \frac{\alpha_1 -1}{2} y \sqrt{\frac{2}{m(1+\alpha_1)}} \right)^{\frac{2}{1-\alpha_1}}.
\end{array}
\right.
$$

\section{Singular Limits of the Interface}
For any small $\varepsilon>0$, denote the solution of \eqref{P} with initial data in \eqref{ini-1} and
\eqref{ini-2} by $u^\varepsilon (x,t)$. Now on the time interval $(0,T)$ for some $T>0$, we consider the limit of the $a$-level set of $u^\varepsilon$. Using the maximum principle we have $0\leq u^\varepsilon(x,t)<1$. We will use the following notation.
$$
\begin{array}{c}
\Omega^\varepsilon_+ (t) := \{x\mid u^\varepsilon (x,t)>a\},\qquad  \Omega^\varepsilon_- (t) := \{x\mid u^\varepsilon (x,t)<a\},\\
\Gamma^\varepsilon (t) := \{x\mid u^\varepsilon(x,t)=a\},\qquad 
d^\varepsilon (x,t):= (-1)^\mu \mbox{dist} (x,\Gamma^\varepsilon(t)),
\end{array}
$$
where $\mu=1$ if $x\in \Omega^\varepsilon_- (t)$ and $\mu=2$ if $x\in \Omega^\varepsilon_+ (t)$.
So $d^\varepsilon (x,t)$ is a signed distance function with respect to $\Gamma^\varepsilon(t)$, and when $\Gamma^\varepsilon(t)$ is a smooth hypersurface,  
\begin{equation}
\left. \nabla d^\varepsilon (x,t)\right|_{x\in \Gamma^\varepsilon(t)} = {\bf n}(x)|_{x\in \Gamma^\varepsilon (t)},\qquad
\left. \Delta d^\varepsilon (x,t)\right|_{x\in \Gamma^\varepsilon(t)} = -(N-1) \kappa(x)|_{x\in \Gamma^\varepsilon (t)},\qquad
\end{equation}
where ${\bf n}(x)$ denotes the unit normal vector at $x\in \Gamma^\varepsilon(t)$ pointing to $\Omega^\varepsilon_+(t)$ from $\Omega^\varepsilon_-(t)$, $\kappa=\kappa(x)$ denotes the mean curvature of
$\Gamma^\varepsilon(t)$ at $x\in \Gamma^\varepsilon(t)$, which is positive if $\Omega^\varepsilon_+(t)$ is convex near $x$.
By the standard theory for porous media equation (\cite{Vaz-book, Wu-book}) we see that, $u^\varepsilon(x,t)$ is continuous in $\overline{\Omega}\times [0,T)$ and classical in the domain where $u^\varepsilon(x,t)>0$. Hence, $\Gamma^\varepsilon(t)$ consists of $(N-1)$-dimensional hypersurfaces, which are called the interfaces of $u^\varepsilon(\cdot,t)$.
As in \cite{NMHS, Lou2007}, we use the matched asymptotic expansion to derive the motion law of the interfaces as $\varepsilon\to 0$.

In what follows, we assume $f\in C^2([0,1])$ and that it satisfies the additional condition:
$$
f'(0)<0 \mbox{ or } f''(0)\not= 0,\qquad f'(1)<0 \mbox{ or } f''(1)\not= 0.
$$

\subsection{The outer expansion}
If $x$ is far (comparing with $\varepsilon$) from the interface, the solution is considered as a normal one which has no sharp layer. So, we assume
$$
u^\varepsilon = u_0 + \varepsilon u_1 + \varepsilon^2 u_2 +\cdots.
$$
Substituting it into the equation and collecting the terms of $O(\varepsilon^{-2}),\ O(\varepsilon^{-1})$ and $O(1)$ we have
\begin{equation}\label{3-eqs}
f(u_0)=0,\qquad f'(u_0)u_1=0,\qquad  u_{0t}=\Delta u_0^m + {\bf q}(x,t) \nabla u_0^k + f'(u_0) u_2 +\frac{f''(u_0)}{2} u_1^2.
\end{equation}
Since we are considering $x$ which is far from $\Gamma^\varepsilon (t)$, we see from the first equality that
$$
u_0 (x,t)=0 \mbox{ when } x\in \Omega^\varepsilon_-(t),\qquad
u_0 (x,t)=1 \mbox{ when } x\in \Omega^\varepsilon_+(t).
$$
In $\Omega^\varepsilon_-(t)$, if $f'(0)<0$ then $u_1=0$, if $f'(0)=0$ but $f''(0)\not= 0$, then by the third equation in \eqref{3-eqs} we also have $u_1=0$. Similarly, we have $u_1 =0$ in $\Omega^\varepsilon_+(t)$.

\subsection{The inner expansion}
Now we consider the points near $\Gamma(t)$. To specify the behavior clearer we use the rescaled variable
$$
\xi := \frac{d^\varepsilon (x,t)}{\varepsilon},
$$
and assume 
\begin{equation}\label{inner-exp}
d^\varepsilon =d_0 +\varepsilon d_1 +\varepsilon^2 d_2 +\cdots,\qquad 
u^\varepsilon = U_0(x,\xi,t)+\varepsilon U_1(x,\xi,t)+\varepsilon^2 U_2(x,\xi,t)+\cdots.
\end{equation}
Substituting them into the equation we see that the coefficient of $\varepsilon^{-2}$ is
\begin{equation}\label{inner-1}
(U_0^m)_{\xi\xi} + f(U_0)=0,
\end{equation}
the coefficient of $\varepsilon^{-1}$ is 
\begin{equation}\label{inner-2}
 m \big( U_0^{m-1}U_1\big)_{\xi\xi} + f'(U_0)U_1 = U_{0\xi} d_{0t} - (U_0^m)_{\xi} \Delta d_0 - {\bf q}\cdot  \nabla d_0 \big(U_0^k\big)_{\xi}.
\end{equation}
In order to make the inner expansion match the outer one, we require the matching conditions:
\begin{equation}\label{match-cond}
U_0(x,-\infty,t)=0,\quad U_0(x,\infty,t)=1,\quad U_1(x,\pm \infty,t)=0,
\end{equation}
as well as the normalization conditions:
\begin{equation}\label{normal-cond}
U_0(x,0,t)=a,\quad U_1(x,0,t)=0.
\end{equation}

Now we solve the solutions $U_0$ and $U_1$. Clearly, the solution $\Phi$ of \eqref{ss}-\eqref{ini-p} solves \eqref{inner-1}-\eqref{match-cond}-\eqref{normal-cond}. So, for any given $(x,t)$,
$$
U_0(x,\xi,t)\equiv \Phi(\xi).
$$
Substituting this result into \eqref{inner-2} we have
\begin{equation}\label{inner-3}
 m \big( \Phi^{m-1}U_1\big)_{\xi\xi} + f'(\Phi)U_1 = \Phi_{\xi} d_{0t} - (\Phi^m)_{\xi} \Delta d_0 - {\bf q}\cdot  \nabla d_0 \big(\Phi^k\big)_{\xi}.
\end{equation}
Note that we assume $f\in C^2([0,1])$ in this part, and so $f(u)\sim -u^{\alpha_0}$ for some $\alpha_0\geq 1$. Hence \eqref{varphi-sim-eq} and \eqref{decay-0} implies that, when $0\leq \alpha_0 <m$, as $y\to y_*>-\infty$, there holds 
$$
\Phi^{2m-3}\Phi^2_{\xi}\to 0, \qquad \Phi^{2m-2}\Phi_{\xi}\to 0, \qquad \Phi^{2m-2}\Phi_{\xi\xi}\to 0.
$$
Therefore,   
\begin{equation}\label{part-1}
m\int_{y_*}^\infty [\Phi^{m-1}\Phi_\xi] \big(\Phi^{m-1}U_1\big)_{\xi\xi} = m \int_{y_*}^\infty \big(\Phi^{m-1}\Phi_\xi \big)_{\xi\xi} \Phi^{m-1} U_1.
\end{equation}
Multiplying the equation \eqref{inner-3} with $\Phi^{m-1}\Phi_\xi$, integrating it over $\xi\in [y_*,\infty)$ and using \eqref{part-1} we have
\begin{equation}\label{part-2}
\int_{y_*}^{\infty} U_1 \Phi^{m-1} \Big[ 
\big(m \Phi^{m-1}\Phi_\xi\big)_{\xi\xi}  + \Phi_{\xi} f'(\Phi) \Big] d\xi
= A d_{0t} - B \Delta d_0 - C {\bf q}\nabla d_0, 
\end{equation}
where 
$$
A:= \int_{y_*}^\infty \Phi^{m-1} \Phi_{\xi}^2 d\xi,\qquad 
B:= \int_{y_*}^\infty (\Phi^m)_\xi \Phi^{m-1} \Phi_{\xi} d\xi,\qquad
C:= \int_{y_*}^\infty (\Phi^k)_\xi \Phi^{m-1} \Phi_{\xi} d\xi
$$
are all positive constants depending only on $f$ and $m$. Differentiating $(\Phi^m)_{\xi\xi}+f(\Phi)=0$ we see that the left hand side of \eqref{part-2} is zero, and so we have
$$
Ad_{0t} = B\Delta d_0 + C{\bf q} \nabla d_0.
$$
This gives the motion law of the interface as $\varepsilon\to 0$:
\begin{equation}\label{MCF-eq}
V = d_{0t} = \frac{B}{A}(N-1)\kappa (x,t) + \frac{C}{A}{\bf q}(x,t)\cdot {\bf n},\qquad x\in \Gamma(t),
\end{equation}
where $\Gamma(t)$ is the limiting interface of $\Gamma^\varepsilon(t)$ as $\varepsilon\to 0$, and $V$  denotes the normal velocity of $\Gamma(t)$ along ${\bf n}$.

\section{The Boundary Conditions}
Denote by $\Lambda(t)$ the contact $(N-2)$-dimensional hypersurface between $\Gamma(t)$ and $\partial \Omega$. Now we use the inner expansion to derive the boundary condition on $\Lambda(t)$.

Taking $\tilde{x}\in \Gamma^\varepsilon(t)$ sufficiently close to $\Lambda(t)$, and substituting \eqref{inner-exp} into the boundary condition of $u^\varepsilon$ in \eqref{P} we have
$$
\frac{1}{\varepsilon} \Big( (U_0^m)_\xi \nabla d_0\Big) \cdot \nu + \nabla (U_0^m) + m\Big(U_0^{m-1}U_1\Big)_{\xi} \nabla d_0 +\cdots = \frac{1}{\varepsilon} G(\tilde{x},t, U_0+\varepsilon U_1 +\cdots).
$$
So, the coefficients of $\varepsilon^{-1}$ satisfy
$$
 (\Phi^m)_\xi {\bf n}\cdot \nu = G(\tilde{x},t, \Phi).
$$
Multiplying it with $\Phi_\xi$ and integrating over $\xi \in [y_*, \infty)$ we have 
$$
D {\bf n}\cdot \nu = g(\tilde{x},t) := \int_0^1 G(\tilde{x},t,r)dr,
$$
with 
$$
D:= \int_{y_*}^{\infty} (\Phi^m)_{\xi} \Phi_\xi d\xi.
$$
This gives the boundary condition for the moving interface $\Gamma(t)$ on the contacting hypersurface $\Lambda(t)$:
\begin{equation}\label{limit-BC}
{\bf n}\cdot \nu = \frac{1}{D}g(x,t),\qquad x \in \Lambda(t).
\end{equation}

As a conclusion, we can see that the singular limit as $\varepsilon\to 0$ of the boundary value problem  \eqref{P} is a mean curvature flow $\Gamma(t)$, which moves according to \eqref{MCF-eq} with nonlinear boundary condition \eqref{limit-BC}. 

In the special case where $\Omega(t)=\omega \times \R$ for some bounded domain $\omega\subset \R^{N-1}$ with smooth boundary, and $\Gamma(t)$ is a graphic hypersurface, that is, the graph of a function $x_N = w(x',t)\ (x'\in \omega)$. If we assume $\Omega^\varepsilon_+$ is the upper half of the cylinder $\Omega$, then
$$
{\bf n} = \frac{(-Dw,1)}{\sqrt{1+|Dw|^2}},\qquad 
(N-1)\kappa = \mbox{div} {\bf n} = \left( \delta_{ij}-\frac{D_i w D_j w}{1+|Dw|^2}\right) \frac{D_{ij}w}{\sqrt{1+|Dw|^2}},
$$
where $D_i w= \frac{\partial w}{\partial x_i}\ (i=1,2,\cdots, N-1)$, $Dw =(D_1 w,\cdots, D_{N-1}w)$. Hence,   \eqref{MCF-eq}-\eqref{limit-BC} is converted to
\begin{equation}\label{graph-MCF}
\left\{
\begin{array}{ll}
\displaystyle w_t  =  \displaystyle \frac{B}{A}\left( \delta_{ij}-\frac{D_i w D_j w}{1+|Dw|^2}\right) D_{ij}w + \frac{C}{A} {\bf q}(-Dw,1),  & x'\in \omega,\ t>0,\\
\displaystyle \frac{-Dw \cdot \nu' }{\sqrt{1+|Dw|^2}} = g(x',w(x',t),t), &  x'\in \partial \omega,\ t>0,
\end{array}
\right.
\end{equation}
where $\nu'$ is the outward unit normal vector of $\partial \omega$, and so $\nu = (\nu',0)$ in such special case.

\end{document}